\documentclass[11pt,reqno]{amsart}
\usepackage{amssymb,amscd,amsbsy}
\usepackage{amssymb,amscd,amsbsy,mathrsfs}
\newcommand{\lb}{\linebreak}

\newcommand{\z}{\zeta}

\renewcommand{\l}{\lambda}

\newcommand{\s}{\sigma}

\newcommand{\f}{\varphi}
\renewcommand{\o}{\omega}

\newcommand{\D}{\Delta}
\renewcommand{\L}{\Lambda}

\newcommand{\E}{{\mathscr E}}
\newcommand{\cd}{{\mathscr D}}

\newcommand{\h}{{\mathscr H}}

\newcommand{\K}{{\mathscr K}}

\newcommand{\X}{{\mathscr X}}
\newcommand{\Y}{{\mathscr Y}}

\newcommand{\C}{{\Bbb C}}
\newcommand{\T}{{\Bbb T}}

\newcommand{\dd}{{\Bbb D}}
\newcommand{\R}{{\Bbb R}}
\newcommand{\Z}{{\Bbb Z}}

\newcommand{\bS}{{\boldsymbol S}}

\newcommand{\rf}[1]{(\ref{#1})}

\newcommand{\df}{\stackrel{\mathrm{def}}{=}}

\newcommand{\spn}{\operatorname{span}}

\newcommand{\clos}{\operatorname{clos}}

\newcommand{\const}{\operatorname{const}}

\newcommand{\eeq}{\end{equation}}
\newcommand{\beq}{\begin{equation}}
\newcommand{\bay}{\begin{eqnarray}}
\newcommand{\ba}{\begin{align*}}
\newcommand{\ea}{\end{align*}}
\newcommand{\ey}{\end{eqnarray}}
\newcommand{\bey}{\begin{eqnarray*}}
\newcommand{\eey}{\end{eqnarray*}}

\newcommand{\be}{\infty}

\newcommand{\bl}{\blacksquare}

\newcommand{\Range}{\operatorname{Range}}

\newcommand{\Pf}{{\bf Proof. }}
\newcommand{\im}{\operatorname{Im}}

\newtheorem{thm}{\hspace{\parindent}Theorem}[section]

\theoremstyle{remark}

\newtheorem*{rem*}{Remark}

\newcommand{\CA}{{\rm C}_{\rm A}}
\newcommand{\OL}{{\rm OL}}

\newcommand{\fF}{{\frak F}}

\newcommand{\OLA}{{\rm OL}_{\rm A}}
\newcommand{\LIA}{\Li_{\rm A}}
\newcommand{\CAo}{{\rm C}_{{\rm A},1}}
\newcommand{\CAbe}{{\rm C}_{{\rm A},\be}}

\newcommand\Li{{\rm Lip}}

\newcommand\dg{\frak D}

\newcommand{\ri}{{\rm i}}

\begin{document}

\newcommand{\vse}{\vspace{.2in}}
\numberwithin{equation}{section}

\title{Dissipative operators and operator Lipschitz functions}
\author{A.B. Aleksandrov and V.V. Peller}
\thanks{The paper was prepared under the support of the program
of the Presidium of the RAS N\textsuperscript{\underline{o}}01 "Fundamental Mathematics and its Applications"
 (grant PRAS-18-01), the research of the first author is also supposed by RFBR grant 17-01-00607 .
The publication was prepared with the support of the
"RUDN University Program 5-100"}
\thanks{Corresponding author: V.V. Peller; email: peller@math.msu.edu}

\begin{abstract}
The purpose of this paper is to obtain an integral representation for the difference $f(L_1)-f(L_2)$ of functions of maximal dissipative operators. This representation in terms of double operator integrals will allow us to establish Lipschitz type estimates for functions of maximal dissipative operators. We also consider a similar problem for quasicommutators, i.e., operators of the form $f(L_1)R-Rf(L_2)$.
\end{abstract}

\maketitle

%\footnotesize
%\tableofcontents

\

\setcounter{section}{0}
\section{\bf Introduction}
\setcounter{equation}{0}
\label{In}

\

It was discovered in Farforovskaya's \cite{F} paper  that a Lipschitz function $f$ on the real line $\R$ (i.e., $f$ satisfies the inequality 
$|f(s)-f(t)|\le\const|s-t|$, $s,\,t\in\R$) does not have to satisfy the inequality
\bay
\label{OpLiAB}
\|f(A)-f(B)\|\le\const\|A-B\|
\ey
for arbitrary bounded self-adjoint operators $A$ and $B$. Functions satisfying 
\rf{OpLiAB} are called {\it operator Lipschitz functions} on $\R$. We refer the reader to \cite{AP2} for a recent detailed survey of operator Lipschitz functions.

Operator Lipschitz functions can be defined on arbitrary closed subset $\fF$ of the complex plane. A complex-valued function $f$ on $\fF$ is called {\it operator Lipschitz} if
\bay
\label{OpLiN1N2}
\|f(N_1)-f(N_2)\|\le\const\|N_1-N_2\|
\ey
for arbitrary bounded normal operators $N_1$ and $N_2$ whose spectra are contained in
$\fF$. It turns out that if $f$ is an operator Lipschitz function on 
a closed unbounded set $\fF$, then
inequality \rf{OpLiN1N2} also holds for not necessarily bounded normal operators $N_1$ and $N_2$ with spectra in $\fF$, see Theorem 3.2.1 of \cite{AP2}.

In this paper we consider the class of operator Lipschitz functions on the closed upper half-plane $\clos\C_+\df\{\z\in\C:~\im\z\ge0\}$ and among operator Lipschitz functions on $\clos\C_+$ we consider those that are analytic in the open half-plane $\C_+$.
We denote the class of such functions by $\OLA(\C_+)$.

The main purpose of this paper is to show that the class $\OLA(\C_+)$ can be characterized as the maximal class of functions $f$, for which the Lipschitz-type estimate
\bay
\label{otsetiLi}
\|f(L_1)-f(L_2)\|\le\const\|L_1-L_2\|
\ey
holds for arbitrary maximal dissipative operators $L_1$ and $L_2$ with bounded difference.

If $f$ is an operator Lipschitz function on the real line and $A$ and $B$ be self-adjoint operators with bounded $A-B$, then the following formula holds
\bay
\label{dvoopin}
f(A)-f(B)=\iint_{\R\times\R}(\dg f)(x,y)\,dE_A(x)(A-B)\,dE_B(y),
\ey
where $E_A$ and $E_B$ are the spectral measures of $A$ and $B$ and the divided difference $\dg f$ is defined on $\R\times\R$ by
\bay
\label{razrakrsamo}
(\dg f)(x,y)\df
\left\{\begin{array}{ll}\displaystyle{\frac{f(x)-f(y)}{x-y}},&x\ne y,\\[.4cm]
f'(x),&x=y,
\end{array}\right.
\ey
(see \cite{BS1} and the survey \cite{AP2}). Recall that by a result of Johnson and Williams \cite{JW}, the operator Lipschitz functions are differentiable everywhere on $\R$ (see also the survey \cite{AP2}, \S\:3.3). The expression on the right-hand side
of \rf{dvoopin} is a double operator integral, see \S\:\ref{dvjoinnLfu}.

In \S\:\ref{glarez} we obtain an analog of this results for maximal dissipative operators. This will allow us to prove a Lipschitz type estimate analogous to \rf{otsetiLi}
for maximal dissipative operators $L_1$ and $L_2$ and functions $f$ of class $\OL_A(\C_+)$.

We obtain in \S\:\ref{komrez} similar results for quasicommutators
$f(L_1)R-Rf(L_2)$ for maximal dissipative operators $L_1$ and $L_2$ and bounded operators $R$.

We give a brief introduction to dissipative operators in \S\:\ref{dissi} and to double operator integrals in \S\:\ref{dvjoinnLfu}.

Note that the proofs of such results for dissipative operators are considerably more complicated than in the case of self-adjoint operators, unitary operators or contractions. First of all, we deal with unbounded functions of unbounded operators. Secondly, unlike in the case of self-adjoint operators, we cannot use spectral projections onto subspaces on which the operators are bounded.

\

\section{\bf Dissipative operators}
\setcounter{equation}{0}
\label{dissi}

\

In this section we give a brief introduction in dissipative operators. We refer the reader to \cite{SNF}, \cite{So} and \cite{AP1} for more detailed information.

Let $\h$ be a Hilbert space. Recall that an operator $L$ (not necessarily bounded)
with dense domain $\cd_L$ in $\h$ is called {\it dissipative} if
$$
\im(Lu,u)\ge0,\quad u\in\cd_L.
$$
A dissipative operator is called {\it maximal dissipative} if it has no proper dissipative extension.

The {\it Cayley transform} of a dissipative operator $L$ is defined by
$$
T\df(L-{\rm i}I)(L+{\rm i}I)^{-1}
$$
with domain $\cd_T=(L+{\rm i}I)\cd_L$ and range $\Range T=(L-{\rm i}I)\cd_L$
(the operator $T$ is not densely defined in general). It can easily be shown that $T$ is a contraction, i.e., $\|Tu\|\le\|u\|$, $u\in\cd_T$, $1$ is not an eigenvalue of $T$, and $\Range(I-T)\df\{u-Tu:~u\in\cd_T\}$ is dense.

Conversely, if $T$ is a contraction defined on its domain $\cd_T$, $1$ is not an eigenvalue of $T$, and $\Range(I-T)$ is dense, then it is the Cayley transform of a dissipative operator $L$ and $L$ is the inverse Cayley transform of $T$:
$$
L={\rm i}(I+T)(I-T)^{-1},\quad \cd_L=\Range(I-T).
$$

A dissipative operator is maximal if and only if the domain of its Cayley transform is the whole Hilbert space.

Every dissipative operator has a maximal dissipative extension.
Every maximal dissipative operator is necessarily closed.

%If $L$ is a maximal dissipative operator, then $-L^*$ is also maximal dissipative.
%Moreover, it is maximal dissipative.

If $L$ is a maximal dissipative operator, then its spectrum $\s(L)$ is contained in the closed upper half-plane $\clos\C_+$.
%and
%\bay
%\label{res}
%\big\|(L-\l I)^{-1}\big\|\le\frac1{|\im\l|},\quad\im\l<0.
%\ey

If $L_1$ and $L_2$ are maximal dissipative operators, we say that the {\it operator $L_1-L_2$ is bounded} if there exists a bounded operator $K$ such that $L_2=L_1+K$.

We proceed now to the construction of functional calculus for dissipative operators. Let $L$ be a maximal dissipative operator and let $T$ be its Cayley transform. Consider its minimal unitary dilation $U$, i.e., $U$ is a unitary operator defined on a Hilbert space $\K$ that contains $\h$ such that
$$
T^n=P_\h U^n\big|\h,\quad n\ge0,
$$
and $\K=\clos\spn\{U^nh:~h\in\h,~n\in\Z\}$. Since $1$ is not an eigenvalue of $T$, it follows that $1$ is not an eigenvalue of $U$ (see \cite{SNF}, Ch. II, \S\,6).

The Sz.-Nagy--Foia\c s functional calculus allows us to define a functional calculus for $T$ on the Banach algebra
$$
\CAo\df\big\{g\in H^\be:~g\quad
\mbox{is continuous on}\quad\T\setminus\{1\}~\big\}.
$$
If $g\in\CAo$, we put
$$
g(T)\df P_\h g(U)\Big|\h.
$$
This functional calculus is linear and multiplicative and
$$
\|g(T)\|\le\|g\|_{H^\be},\quad g\in\CAo,
$$
(see \cite{SNF}, Ch. III).

We can define now a functional calculus for our dissipative operator on the
Banach algebra
$$
\CAbe=\big\{f\in H^\be(\C_+):~f\quad\mbox{is continuous on}\quad
\R\big\}.
$$
Indeed, if $f\in\CAbe$, we put
$$
f(L)\df \big(f\circ\o\big)(T),
$$
where $\o$ is the conformal map of $\dd$ onto $\C_+$ defined by
$\o(\z)\df{\rm i}(1+\z)(1-\z)^{-1}$, $\z\in\dd$.

We can extend now this functional calculus to the class $\LIA(\C_+)$ of Lipschitz functions on $\clos\C_+$ that are analytic in $\C_+$. Suppose that $f\in\LIA(\C_+)$.
We have
\bay
\label{fri}
f(\z)=\frac{f_\ri(\z)}{(\z+\ri)^{-1}},\quad\z\in\C_+,
\qquad
\mbox{where}\quad f_\ri(\z)\df\frac{f(\z)}{\z+\ri}.
\ey
Since $f\in\LIA(\C_+)$, the function $f_\ri$ belongs to $\CAbe$. The (possibly unbounded) operator $f(L)$ can be defined by
\bay
\label{f(L)}
f(L)\df(L+\ri I)f_\ri(L)
\ey
(see \cite{SNF}, Ch. IV, \S\;1). It follows from Th. 1.1 of Ch. IV of \cite{SNF} that 
\bay
\label{fi(L)}
f(L)\supset f_\ri(L)(L+\ri I),
\ey
and so $D(f(L))\supset D(L)$.

We proceed now to the definition of a resolvent self-adjoint dilation of a maximal dissipative operator. If $L$ is a maximal dissipative operator on a Hilbert space $\h$, we say that a self-adjoint operator $A$ in a Hilbert space
$\K$, $\K\supset\h$, is called a {\it resolvent self-adjoint dilation} of $L$
if
$$
(L-\l I)^{-1}=P_\h(A-\l I)^{-1}\Big|\h,\quad \im\l<0.
$$
The dilation is called {\it minimal} if
$$
\K=\clos\spn\big\{(A-\l I)^{-1}v:~v\in\h,~\im\l<0\big\}.
$$
If $f\in C_{A,\be}$, then
$$
f(L)=P_\h f(A)\Big|\h,\quad f\in C_{A,\be}.
$$

A minimal resolvent self-adjoint dilation of a maximal dissipative operator always exists (and is unique up to a natural isomorphism). Indeed, it suffices to take a minimal unitary dilation of the Cayley transform of this operator and
apply the inverse Cayley transform to it.

%Sometimes mathematicians use the term "self-adjoint dilation" rather than
%"resolvent self-adjoint dilation". However, we believe that the term "self-adjoint dilation" is misleading.

Let us define now the semi-spectral measure of a maximal dissipative operator $L$. Let $T$ be the Cayley transform of $L$ and let $\E_T$ be the semi-spectral measure of $T$ on the unit circle $\T$ defined by
$$
\E_T(\D)\df P_\h E_U(\D),
$$
where $\D$ is a Borel subset of $\T$ and $E_U$ is the spectral measure of the minimal unitary dilation $U$ of $T$.

Then
\bay
\label{ssu}
g(T)=\int_\T g(\z)\,d\E_T(\z),\quad g\in C_{A,1}.
\ey
We can define now the semi-spectral measure $\E_L$ of $L$ by
$$
\E_L(\D)=\E_T\big(\o^{-1}(\D)\big),\quad\D\quad\mbox{is a Borel subset of}
\quad\R.
$$
It follows easily from \rf{ssu} that
\bay
\label{fL}
f(L)=\int_\R f(x)\,d\E_L(x),\quad f\in C_{A,\be}.
\ey
It can easily be verified that if $A$ is the minimal self-adjoint resolvent dilation of $L$ and $E_A$ is the spectral measure of $A$, then
$$
\E_L(\D)=P_\h E_A(\D)\big|\h,\quad\D\quad\mbox{is a Borel subset of}
\quad\R.
$$

\

\section{\bf Double operator integrals and operator Lipschitz functions}
\setcounter{equation}{0}
\label{dvjoinnLfu}

\

Double operator integrals
$$
\iint_{\X\times\Y}\Phi(x,y)\,dE_1(x)Q\,dE_2(y)
$$
were introduced by Yu.L. Daletskii and S.G. Krein, see \cite{DK}. Later
Birman and Solomyak elaborated their beautiful theory of double operator integrals \cite{BS1} and \cite{BS2} (see also \cite{AP2} and references therein). Here $\Phi$ is a bounded measurable function, $E_1$ and $E_2$ are spectral measures on a separable Hilbert space
defined on $\s$-algebras of subsets of sets $\X$ and $\Y$ and
$Q$ is a bounded linear operator. 

The approach of Birman and Solomyak  \cite{BS1} starts with the case when $Q\in\bS_2$, i.e., $Q$ is a Hilbert--Schmidt operator.
Under this assumption double operator integrals can be defined for arbitrary bounded measurable functions $\Phi$. Put
$$
{\mathcal E}(\L\times\D)Q=E_1(\L)QE_2(\D),\quad Q\in\bS_2,
$$
where $\L$ and $\D$ are measurable subsets of $\X$ and $\Y$. Clearly,
${\mathcal E}$ takes values in the set of orthogonal projections on the Hilbert space $\bS_2$. 
It was shown in \cite{BS1} (see also \cite{BS3}) that ${\mathcal E}$ extends to a spectral measure on
$\X\times\Y$. If $\Phi$ is a bounded measurable function on $\X\times\Y$, then
$$
\iint_{\X\times\Y}\Phi(x,y)\,d E_1(x)Q\,dE_2(y)\df
\left(\,\iint_{\X\times\Y}\Phi\,d{\mathcal E}\right)Q.
$$
Clearly,
$$
\left\|\iint_{\X\times\Y}\Phi(x,y)\,dE_1(x)Q\,dE_2(y)\right\|_{\bS_2}
\le\|\Phi\|_{L^\be}\|Q\|_{\bS_2}.
$$

If $Q$ is an arbitrary bounded operator, then for the double operator integral to make sense, $\Phi$ has to be a Schur multiplier with respect to $E_1$ and $E_2$, (see \cite{Pe1} and \cite{AP2}). It is well known (see \cite{Pe1}, \cite{AP2} and \cite{Pi}) that $\Phi$ is a Schur multiplier if and only if 
$\Phi$ belongs to the Haagerup tensor product 
$L^\be_{E_1}\otimes_{\rm h}L^\be_{E_2}$, i.e.,
admits a representation
\bay
\label{htenppre}
\Phi(x,y)=\sum_{n\ge0}\f_n(x)\psi_n(y),
\ey
where the $\f_n$ and $\psi_n$ are measurable functions such that
$$
\sum_{n\ge0}|\f_n|^2\in L^\be_{E_1}\quad\mbox{and}\quad
\sum_{n\ge0}|\psi_n|^2\in L^\be_{E_2}.
$$
In this case
$$
\iint\Phi(x,y)\,dE_1(x)Q\,dE_2(y)=
\sum_{n\ge0}\left(\int\f_n(x)\,dE_1(x)\right)Q\left(\int\psi_n(y)\,dE_2(y)\right)
$$
and the right-hand side does not depend on the choice of a representation in 
\rf{htenppre}. 

In this paper we need double operator integrals with respect to {\it semi-spectral measures}
\bay
\label{dvoopipolu}
\iint\Phi(x,y)\,d\E_1(x)Q\,d\E_2(y).
\ey
Such double operator integrals were introduced in \cite{Pe2} (see also \cite{Pe3}).
By analogy with the case of double operator integrals with respect to spectral measures, double operator integrals of the form \rf{dvoopipolu} can be defined for arbitrary bounded measurable functions $\Phi$ in the case when $Q\in\bS_2$ and for 
functions $\Phi$ in $L^\be_{\E_1}\otimes_{\rm h}L^\be_{\E_2}$ in the case of an arbitrary bounded operator $Q$.

We refer the reader to the recent surveys \cite{AP2} and \cite{Pe4} for detailed information.

Suppose now that $f\in\OLA(\C_+)$. It follows from a result of Johnson and Williams
\cite{JW} (see also the survey \cite{AP2}, \S\;3.3) that $f$ is differentiable everywhere on $\clos\C_+$. We define the divided difference $\dg f$ of $f$ on
$\clos\C_+\times\clos\C_+$ by
\bay
\label{razrakr}
(\dg f)(z,w)\df
\left\{\begin{array}{ll}\displaystyle{\frac{f(z)-f(w)}{z-w}},&z\ne w,\\[.4cm]
f'(z),&z=w.
\end{array}\right.
\ey

Let $\CA(\C_+)$ be the set of functions analytic in $\C_+$ and continuous in $\clos\C_+$ and having a finite limit at infinity.
We need the following characterization of the divided differences 
$\dg f$ for functions in $\OLA(\C_+)$ (see \cite{AP2}, Theorems 3.9.6):

\medskip

{\it Let $f$ be a function continuous on $\clos\C_+$ and analytic in $\C_+$. Then $f\in\OLA(\C_+)$ if and only if
$\dg f$ admits a representation 
\bay
\label{Haagerpre}
(\dg f)(z,w)=\sum_{n\ge1}\f_n(z)\psi_n(w),\quad z,\;w\in\clos\C_+,
\ey
where $\f_n, \psi_n\in \CA(\C_+)$   such that
\bay
\label{proisumm}
\left(\sup_{z\in\C_+}\sum_{n\ge1}|\f_n(z)|^2\right)^{1/2}\left(\sup_{w\in\C_+}\sum_{n\ge1}|\psi_n(w)|^2\right)^{1/2}<\be.
\ey
If $f\in\OLA(\C_+)$, then the functions $\f_n$ and $\psi_n$ can be chosen so that
the left-hand side of {\em\rf{proisumm}} is equal to $\|f\|_{\OLA}$.}

Recall that
$$
\|f\|_{\OLA}\df\sup\left\{\frac{\|f(N_1)-f(N_2)\|}{\|N_1-N_2\|}\right\},
$$
the supremum being taken over all distinct normal operators $N_1$ and $N_2$ with spectra in $\clos\C_+$ such that the operator $N_1-N_2$ is bounded.

\medskip

This allows us to consider double operator integrals
$$
\iint_{\R\times\R}(\dg f)(x,y)\,d\E_1(x)Q\,d\E_2(y),
$$
where $\E_1$ and $\E_2$ are semi-spectral measures and $Q$ is a bounded linear operator, and this double operator integral is equal to
$$
\sum_{n\ge1}\left(\int_\R\f_n\,d\E_1\right)Q\left(\int_\R\psi_n\,d\E_2\right),
$$
where the functions $\f_n$ and $\psi_n$ satisfy \rf{Haagerpre} and \rf{proisumm}.
Moreover, the following inequality holds:
\bay
\label{otsdvopi}
\left\|\iint_{\R\times\R}(\dg f)(x,y)\,d\E_1(x)Q\,d\E_2(y)\right\|
\le\|f\|_{\OLA}\|Q\|.
\ey

Note that it can be deduced from Theorem 2.2.3 in \cite{AP2} that 
if $f\in\LIA(\C_+)$, then $f\in\OLA(\C_+)$ if and only if $f$ is an operator Lipschitz functions on $\R$. 
Moreover, $\|f\|_{\OLA(\C_+)}=\|f\|_{\OL(\R)}$, where
$$
\|f\|_{\OL(\R)}\df\sup\frac{\|f(A)-f(B)\|}{\|A-B\|},
$$
the supremum is being taken over all distinct bounded self-adjoint operators $A$ and $B$.

%In a similar way we can consider the case of functions analytic in the unit disk.
%Let $f\in\OLA(\dd)$, i.e., $f$ is an operator Lipschitz function 
%on $\clos\dd$ and analytic in $\dd$. Then by Theorem 3.9.3 of \cite{AP2}
%there are functions $\f_n$ and $\psi_n$ in the disk-algebra $\CA$
%(i.e., they are continuous on $\clos\dd$ and analytic in $\dd$) such that
%$$
%(\dg f)(z,w)=\sum_{n\ge1}\f_n(z)\psi_n(w),\quad z,\;w\in\clos\dd,
%$$
%and
%$$
%\left(\sup_{z\in\dd}\,\sum_{n\ge1}|\f_n(z)|^2\right)\left(\sup_{w\in\dd}\,\sum_{n\ge1}|\psi_n(w)|^2\right)<\be.
%$$
%Here the divided difference $\dg f$ on $\clos\dd\times\clos\dd$ can be defined by analogy with \rf{razrakr}.
%
%Suppose now that $T_1$ and $T_2$ are contractions on Hilbert space. Then it follows from Theorem 3.9.9 of \cite{AP2}, that
%\bay
%\label{fT1T2sum}
%f(T_1)-f(T_2)=\sum_{n\ge1}\f_n(T_1)(T_1-T_2)\psi_n(T_2)
%\ey
%and the series on the right converges in the weak operator topology.
%

\

\section{\bf The integral representation and Lipschitz type estimates}
\setcounter{equation}{0}
\label{glarez}

\

In this section for function $f$ of class $\OLA(\C_+)$ and for maximal dissipative operators $L_1$ and $L_2$ with bounded difference, we obtain a representation of the difference \lb$f(L_1)-f(L_2)$ in terms 
of a double operator integral. This allows us to obtain a Lipschitz type estimate for functions of maximal dissipative operators.

Let $f\in\OLA(\C_+)$ and let $L$ be a maximal dissipative operator. It can easily be deduced from  \rf{f(L)} and \rf{fi(L)} that $f(L)(L+\ri I)^{-1}=f_\ri(L)$,
where $f_\ri$ is defined in \rf{fri}.

The following theorem holds for arbitrary maximal dissipative operators $L_1$ and $L_2$ without the assumption that $L_1-L_2$ is bounded.

\begin{thm} 
\label{bezep}
Let $L_1$ and $L_2$ be maximal dissipative operators and let $f\in\OLA(\C_+)$.
Suppose that $\f_n$ and $\psi_n$ are functions in $\CA(\C_+)$ 
satisfying {\em\rf{Haagerpre}} and {\rm\rf{proisumm}}. Then
\begin{multline}
\label{dlifo}
f_\ri(L_1)(L_2+\ri I)^{-1}-(L_1+\ri I)^{-1}f_\ri(L_2)\\
=f(L_1)(L_1+\ri I)^{-1}(L_2+\ri I)^{-1}-(L_1+\ri I)^{-1}f(L_2)(L_2+\ri I)^{-1}\\
=\sum_{n=1}^\be\f_n(L_1)\big((L_2+\ri I)^{-1}-(L_1+\ri I)^{-1}\big)\psi_n(L_2)
\end{multline}
and the series on the right converges in the weak operator topology. 
\end{thm} 

\Pf
It follows from \rf{Haagerpre}  that the following identity holds:
\begin{multline}
\label{izw}
f(z)(z+\ri)^{-1}(w+\ri)^{-1}-(z+\ri)^{-1}f(w)(w+\ri)^{-1}\\
=\sum_{n=1}^\be\f_n(z)z(z+\ri)^{-1}(w+\ri)^{-1}\psi_n(w)-
\sum_{n=1}^\be\f_n(z)(z+\ri)^{-1}w(w+\ri)^{-1}\psi_n(w)
\end{multline}
for all $z,\;w\in\clos\C_+$. Let $\E_1$ and $\E_2$ be the semi-spectral measure of $L_1$ and $L_2$. We have
\begin{multline*}
\iint_{\R\times\R}
\big(f(x)(x+\ri)^{-1}(y+\ri)^{-1}-(x+\ri)^{-1}f(y)(y+\ri)^{-1}\big)
\,d\E_1(x)\,d\E_2(y)\\[.2cm]
=f_\ri(L_1)(L_2+\ri I)^{-1}-(L_1+\ri I)^{-1}f_\ri(L_2)\\[.2cm]
=f(L_1)(L_1+\ri I)^{-1}(L_2+\ri I)^{-1}-(L_1+\ri I)^{-1}f(L_2)(L_2+\ri I)^{-1}.
\end{multline*}
Next, 
\begin{align*}
\iint_{\R\times\R}&
\left(\sum_{n=1}^\be\f_n(x)x(x+\ri)^{-1}(y+\ri)^{-1}\psi_n(y)\right)
\,d\E_1(x)\,d\E_2(y)\\[.2cm]
&=\sum_{n=1}^\be\f_n(L_1)L_1(L_1+\ri I)^{-1}(L_2+\ri I)^{-1}\psi_n(L_2)\\[.2cm]
&=\sum_{n=1}^\be\f_n(L_1)\big(I-\ri(L_1+\ri I)^{-1}\big)(L_2+\ri I)^{-1}\psi_n(L_2)
\end{align*}
and
\begin{align*}
\iint_{\R\times\R}&
\left(\sum_{n=1}^\be\f_n(x)(x+\ri)^{-1}y(y+\ri)^{-1}\psi_n(w)\right)
\,d\E_1(x)\,d\E_2(y)\\[.2cm]
&=\sum_{n=1}^\be\f_n(L_1)(L_1+\ri I)^{-1}L_2(L_2+\ri I)^{-1}\psi_n(L_2)\\[.2cm]
&=\sum_{n=1}^\be\f_n(L_1)(L_1+\ri I)^{-1}\big(I-\ri(L_2+\ri I)^{-1}\big)\psi_n(L_2).
\end{align*}
It follows that the integral of the right-hand side of \rf{izw} is equal to
\begin{align*}
\sum_{n=1}^\be
\f_n(L_1)\big(\big(I-\ri(L_1+\ri I)^{-1}\big)(L_2+\ri I)^{-1}-
(L_1+\ri I)^{-1}\big(I-\ri(L_2+\ri I)^{-1}\big)\big)\psi_n(L_2)\\[.2cm]
=\sum_{n=1}^\be
\f_n(L_1)\big((L_2+\ri I)^{-1})-(L_1+\ri I)^{-1}\big)\psi_n(L_2).
\end{align*}
To establish \rf{dlifo}, it remains to equate the integral of the left-hand side of 
\rf{izw} with the integral of the right-hand side. $\bl$

\begin{thm}
\label{intflarafudio}
Let $f\in\OLA(\C_+)$. Then for arbitrary
maximal dissipative operators $L_1$ and $L_2$ with bounded $L_1-L_2$, the following formula holds:
$$
f(L_1)-f(L_2)=\iint_{\R\times\R}(\dg f)(x,y)\,d\E_1(x)(L_1-L_2)\,d\E_2(y),
$$
where $\E_j$ is the semi-spectral measure of $L_j$, $j=1,\,2$.
\end{thm}

To be more precise, we mean that the operator $f(L_1)-f(L_2)$ whose domain contains the dense set $\cd_{L_1}=\cd_{L_2}$  is bounded and extends by continuity to the double operator integral on the right-hand side.

\medskip

\Pf Let $\f_n$ and $\psi_n$, $n\ge1$, be functions in $\CA(\C_+)$
satisfying \rf{Haagerpre} and \rf{proisumm}. Then the double operator integral is well defined and
$$
\iint_{\R\times\R}(\dg f)(x,y)\,d\E_1(x)(L_1-L_2)\,d\E_2(y)=
\sum_{n\ge1}\f_n(L_1)(L_1-L_2)\psi_n(L_2),
$$
where the series converges in the weak operator topology.  Put
$$
Q\df\sum_{n\ge1}\f_n(L_1)(L_1-L_2)\psi_n(L_2).
$$
Then $Q$ is a bounded operator and
\begin{multline*}
(L_1+\ri I)^{-1}Q(L_2+\ri I)^{-1}=\sum_{n\ge1}(L_1+\ri I)^{-1}\f_n(L_1)(L_1-L_2)\psi_n(L_2)(L_2+\ri I)^{-1}\\[.2cm]
=\sum_{n\ge1}\f_n(L_1)(L_1+\ri I)^{-1}\big((L_1+\ri I)-(L_2+\ri I)\big)(L_2+\ri I)^{-1}\psi_n(L_2)\\[.2cm]
=\sum_{n\ge1}\f_n(L_1)\big((L_2+\ri I)^{-1}-(L_1+\ri I)^{-1}\big)\psi_n(L_2)\\[.2cm]
=f(L_1)(L_1+\ri I)^{-1}(L_2+\ri I)^{-1}-(L_1+\ri I)^{-1}f(L_2)(L_2+\ri I)^{-1}.
\end{multline*}

Let us show that
$$
f(L_1)(L_1+\ri I)^{-1}(L_2+\ri I)^{-1}=(L_1+\ri I)^{-1}f(L_1)(L_2+\ri I)^{-1}.
$$
Indeed, $\Range(L_2+\ri I)^{-1}=\cd_{L_2}=\cd_{L_1}$, and so we have to prove that
$$
f(L_1)(L_1+\ri I)^{-1}\big|\cd_{L_1}=
(L_1+\ri I)^{-1}f(L_1)\big|\cd_{L_1}.
$$
This follows easily from the equalities  $f(L_1)(L_1+\ri I)^{-1}=f_\ri(L_1)$
and $f(L_1)=\lb(L_1+\ri I)f_\ri(L_1)$. Thus,
\begin{align*}
(L_1+\ri I)^{-1}Q(L_2+\ri I)^{-1}&=
(L_1+\ri I)^{-1}f(L_1)(L_2+\ri I)^{-1}-(L_1+\ri I)^{-1}f(L_2)(L_2+\ri I)^{-1}\\[.2cm]
&=(L_1+\ri I)^{-1}(f(L_1)-f(L_2))(L_2+\ri I)^{-1}.
\end{align*}
Hence, $Qu=f(L_1)u-f(L_2)u$ for all $u\in \cd_{L_1}=\cd_{L_2}$. $\bl$

\medskip

The following result establishes a Lipschitz type estimate for functions of maximal dissipative operators.

\begin{thm}
Let $f\in\OLA(\C_+)$. Then
$$
\|f(L_1)-f(L_2)\|\le\|f\|_{\OLA}\|L_1-L_2\|
$$
for arbitrary maximal dissipative operators $L_1$ and $L_2$ with bounded $L_1-L_2$.
\end{thm}

\Pf This is an immediate consequence of Theorem \ref{intflarafudio}
and inequality \rf{otsdvopi}. $\bl$

\medskip

{\bf Remark.} In the case when the difference $L_1-L_2$ is a {\it Hilbert--Schmidt operator}, the formula 
$$
f(L_1)-f(L_2)=\iint_{\R\times\R}(\dg f)(x,y)\,d\E_1(x)(L_1-L_2)\,d\E_2(y),
$$
holds {\it for an arbitrary function $f$ in $\LIA(\C_+)$} which leads to the following Lipschitz type estimate in the the Hilbert--Schmidt norm:
$$
\|f(L_1)-f(L_2)\|_{\bS_2}\le\big\{\sup|f'(\z)|:~\im\z>0\big\}\|L_1-L_2\|_{\bS_2}.
$$
This is a consequence of Theorem 6.6 of \cite{AP1} and the obvious equality
$$
f(L)u=\lim_{n\to\be}(f\o_n)(L)u=\lim_{n\to\be}f(L)\o_n(L),\quad u\in\cd_L,
$$
for a maximal dissipative operator $L$, see also \S\:6 of \cite{MNP}. Here 
$$
\o_n(z)\df\frac1{\log n}\log\frac{z+\ri n}{z+\ri},\quad n\ge2,
$$
and $\log$ means the principal branch of logarithm.

\

\section{\bf Commutator Lipschitz estimates}
\setcounter{equation}{0}
\label{komrez}

\

The purpose of this section is to obtain Lipschitz type norm estimates of (quasi) commutators of the form 
$
f(L_1)R-Rf(L_2)
$
in terms of the norms of $L_1R-RL_2$, where $L_1$ and $L_2$ are maximal dissipative operators and $R$ is a bounded linear operator.

We say that the operator $L_1R-RL_2$ is bounded if $R(\cd_{L_2})\subset\cd_{L_1}$ and 
$$
\|L_1Ru-RL_2u\|\le\const\|u\|\quad\text{for every}\quad u\in\cd_{L_2}.
$$

The following result is an analog of Theorem \ref{glarez}. It also holds for arbitrary maximal dissipative operators $L_1$ and $L_2$ and for an arbitrary bounded operator $R$ without any additional assumptions.

\begin{thm} 
Let $L_1$ and $L_2$ be maximal dissipative operators and let $f\in\OLA(\C_+)$ and let $R$ be a bounded operator.
Suppose that $\f_n$ and $\psi_n$ are functions in $\CA(\C_+)$ 
satisfying {\em\rf{Haagerpre}} and {\rm\rf{proisumm}}. Then
\begin{multline}
\label{dlifoR}
f(L_1)(L_1+\ri I)^{-1}R(L_2+\ri I)^{-1}-(L_1+\ri I)^{-1}Rf(L_2)(L_2+\ri I)^{-1}\\
=\sum_{n=1}^\be\f_n(L_1)\big((R(L_2+\ri I)^{-1}-(L_1+\ri I)^{-1}R\big)\psi_n(L_2)
\end{multline}
and the series on the right converges in the weak operator topology. 
\end{thm}

\medskip

\Pf
We are going to use \rf{izw}. We have
\begin{multline*}
\int_\R\int_\R(x+\ri)^{-1}(f(x)-f(y))(y+\ri)^{-1}\,d\E_1(x)R\,d\E_2(y)\\[.2cm]
=f(L_1)(L_1+\ri I)^{-1}R(L_2+\ri I)^{-1}-(L_1+\ri I)^{-1}Rf(L_2)(L_2+\ri I)^{-1}.
\end{multline*}
On the other hand, this is equal to
\begin{multline*}
\int_\R\int_\R\sum_{n=1}^\be
\f_n(x)x(x+\ri)^{-1}(y+\ri)^{-1}\psi_n(y)\,d\E_1(x)R\,d\E_2(y)\\[.2cm]
-\int_\R\int_\R\sum_{n=1}^\be
\sum_{n=1}^\be\f_n(x)(x+\ri)^{-1}y(y+\ri)^{-1}\psi_n(y)\,d\E_1(x)R\,d\E_2(y)\\[.2cm]
=\sum_{n=1}^\be\f_n(L_1)L_1(L_1+\ri I)^{-1}R(L_2+\ri I)^{-1}\psi_n(L_2)\\[.2cm]
-\sum_{n=1}^\be\psi_n(L_1)(L_1+\ri I)^{-1}RL_2(L_2+\ri I)^{-1}\psi_n(L_2)\\[.2cm]
=\sum_{n=1}^\be
\f_n(L_1)\big(I-\ri(L_1+\ri I)^{-1}\big)R(L_2+\ri I)^{-1}\psi_n(L_2)\\[.2cm]
-\sum_{n=1}^\be\psi_n(L_1)(L_1+\ri I)^{-1}R\big(I-\ri(L_2+\ri I)^{-1}\big)\\[.2cm]
=\sum_{n=1}^\be\f_n(L_1)\big(R(L_2+\ri I)^{-1}-(L_1+\ri I)^{-1}R\big)\psi_n(L_2).
\quad\bl
\end{multline*}

\begin{thm}
\label{intprequasiko}
Let $f\in\OLA(\C_+)$. Then for arbitrary bounded operator $R$ and
maximal dissipative operators $L_1$ and $L_2$ with bounded $L_1R-RL_2$, the following formula holds:
$$
f(L_1)R-Rf(L_2)=\iint_{\R\times\R}(\dg f)(x,y)\,d\E_1(x)(L_1R-RL_2)\,d\E_2(y),
$$
where $\E_j$ is the semi-spectral measure of $L_j$, $j=1,\,2$.
\end{thm}

Note that the right-hand side is a bounded operator defined on the whole space, while the left-hand side is a bounded operator defined on $\cd_{L_2}$, which extends by continuity to the operator on the right.

\medskip

\Pf Let $\f_n$ and $\psi_n$, $n\ge1$, be functions in $\CA(\C_+)$
satisfying \rf{Haagerpre} and \rf{proisumm}. Then the double operator integral is well defined and
$$
\iint_{\R\times\R}(\dg f)(x,y)\,d\E_1(x)(L_1R-RL_2)\,d\E_2(y)=
\sum_{n\ge1}\f_n(L_1)(L_1R-RL_2)\psi_n(L_2).
$$
Put
$$
Q\df\sum_{n\ge1}\f_n(L_1)(L_1R-RL_2)\psi_n(L_2).
$$
Then $Q$ is a bounded linear operator and
\begin{multline*}
(L_1+\ri I)^{-1}Q(L_2+\ri I)^{-1}=\sum_{n\ge1}(L_1+\ri I)^{-1}\f_n(L_1)(L_1R-RL_2)\psi_n(L_2)(L_2+\ri I)^{-1}\\[.2cm]
=\sum_{n\ge1}\f_n(L_1)(L_1+\ri I)^{-1}(L_1R-RL_2)(L_2+\ri I)^{-1}\psi_n(L_2)\\[.2cm]
=\sum_{n\ge1}\f_n(L_1)\big((L_1+\ri I)^{-1}L_1R(L_2+\ri I)^{-1}-(L_1+\ri I)^{-1}RL_2(L_2+\ri I)^{-1}\big)\psi_n(L_2).
\end{multline*}
As we have observed in \S\:\ref{glarez}, for a maximal dissipative operator $L$, we have equality $L(L+\ri I)^{-1}=I-\ri(L+\ri I)^{-1}$. Using this equality, we find that
$$
(L_1+\ri I)^{-1}Q(L_2+\ri I)^{-1}
=\sum_{n\ge1}\f_n(L_1)\big(R(L_2+\ri I)^{-1}-(L_1+\ri I)^{-1}R\big)\psi_n(L_2).
$$
By \rf{dlifoR}, the last expression is equal to
\bay
\label{perestavim}
f(L_1)(L_1+\ri I)^{-1}R(L_2+\ri I)^{-1}-(L_1+\ri I)^{-1}Rf(L_2)(L_2+\ri I)^{-1}.
\ey
Now it is time to use the fact that the operator $L_1R-RL_2$ is bounded,
and so \lb$R(\cd_{L_2})\subset\cd_{L_1}$. It follows that 
$\Range R(L_2+\ri I)^{-1}\subset\cd_{L_1}$. Thus, 
$$
f(L_1)(L_1+\ri I)^{-1}R(L_2+\ri I)^{-1}=(L_1+\ri I)^{-1}f(L_1)R(L_2+\ri I)^{-1}.
$$
Hence, the expression in \rf{perestavim} is equal to
\begin{multline*}
(L_1+\ri I)^{-1}f(L_1)R(L_2+\ri I)^{-1}-(L_1+\ri I)^{-1}Rf(L_2)(L_2+\ri I)^{-1}\\
=(L_1+\ri I)^{-1}(f(L_1)R-Rf(L_2))(L_2+\ri I)^{-1}
\end{multline*}
and we can conclude that
$$
(L_1+\ri I)^{-1}Q(L_2+\ri I)^{-1}=(L_1+\ri I)^{-1}(f(L_1)R-Rf(L_2))(L_2+\ri I)^{-1}.
$$
Hence, $Qu=f(L_1)Ru-Rf(L_2)u$ for all $u\in \cd_{L_2}$. $\bl$

\begin{thm}
Let $f\in\OLA(\C_+)$. Then
$$
\|f(L_1)R-Rf(L_2)\|\le\|f\|_{\OLA}\|L_1R-RL_2\|
$$
for an arbitrary bounded operator $R$ and arbitrary maximal dissipative operators $L_1$ and $L_2$ with bounded $L_1R-RL_2$.
\end{thm}

\Pf This is an immediate consequence of Theorem \ref{intprequasiko} and inequality
\rf{otsdvopi}. $\bl$

\

\

%\footnotesize
\noindent
\begin{tabular}{p{7cm}p{15cm}}
A.B. Aleksandrov & V.V. Peller \\
St.Petersburg Branch & Department of Mathematics \\
Steklov Institute of Mathematics  & Michigan State University \\
Fontanka 27, 191023 St.Petersburg & East Lansing, Michigan 48824\\
Russia&USA\\
&and\\
&Peoples' Friendship University\\
& of Russia (RUDN University)\\
&6 Miklukho-Maklaya St., Moscow,\\
& 117198, Russian Federation
\end{tabular}

\end{document}